\documentclass[11pt]{article}

\usepackage{amsthm}

\usepackage[nocompress]{cite}
\usepackage{graphicx}

\usepackage{hyperref}
\hypersetup{
    colorlinks=true,
    linkcolor=blue,
    anchorcolor=blue,
    citecolor=blue,
}

\usepackage{enumitem}

 \usepackage{amsmath}
\usepackage{amsfonts}
\usepackage{amssymb, amsfonts, amsmath, nicefrac}
\usepackage{color}
\usepackage{euscript}

\usepackage[margin=1.0in]{geometry}

\def\sP {{\sf P}}

\newcommand{\cP}{{\mathfrak P}}

\newcommand{\cM}{{\mathfrak M}}

\newcommand{\cF}{{\mathfrak F}}

\newcommand{\cR}{{\mathfrak R}}

\newcommand{\cp}{{\mathfrak p}}

\newtheorem{proposition}{Proposition}[section]

\newtheorem{remark}{Remark}[section]
\newtheorem{definition}{Definition}[section]

\newtheorem{assumption}{Assumption}[section]

\newcommand{\U}{{\cal U}}

\newcommand{\Z}{{\cal Z}}
\newcommand{\F}{{\cal F}}

\newcommand{\A}{{\cal A}}

\newcommand{\X}{{\cal X}}

\newcommand{\s}{{\cal S}}
\newcommand{\R}{{\cal R}}

\newcommand{\be}{\begin{equation}}
\newcommand{\ee}{\end{equation}}

\def\w{\omega}
\def\e{\varepsilon}
\def\O{\Omega}

\newcommand{\avr}{{\sf AV@R}}

\newcommand{\VaR}{{\sf V@R}}

\def\bbr{{\Bbb{R}}} 
\def\bbe{{\Bbb{E}}} 

\def\bbp{{\Bbb{P}}}

\begin{document}

\begin{titlepage}
\title{\bf Statistical robustness and time consistency in risk averse optimization}
\author{{\bf Alexander Shapiro}\thanks{Georgia Institute of Technology, Atlanta, Georgia
30332, USA, \tt{ashapiro@isye.gatech.edu}\newline
Research of this
author was partially supported by the Air Force Office of Scientific Research (AFOSR)
under Grant    FA9550-25-1-0310.}
\and
{\bf Yan Li}\thanks{Texas A\&M University, College Station, Texas 77840, USA, \tt{yan.li@tamu.edu}.}
}


\maketitle
\begin{abstract}
The aim of this manuscript is to investigate risk measures which are not necessarily convex.
We show that law invariant coherent risk measures are not robust with respect to contamination of the probability measure,
in the sense that an arbitrarily small contamination can change their values by an arbitrary amount.
This contrasts with the Value-at-Risk, which is robust against such contamination.
In addition, we show that for a given law invariant risk measure and a given policy,
one can construct, under mild conditions, conditional risk measures with respect to which this policy is time consistent,
in the frameworks of both Markov Decision Processes and Stochastic Optimal Control.
\end{abstract}

\end{titlepage}

\setcounter{equation}{0}
\section{Introduction}
\label{sec-intr}

In this manuscript, we study risk measures which do not necessarily satisfy the convexity axiom,
and their use in sequential decision making under uncertainty.
Modern theory of risk measures was started in \cite{ADEH:1999}, where the axioms of the so-called coherent risk measures were formulated.
An important example of a coherent risk measure is the Average Value-at-Risk,
which can be written as the infimum of a family of expectations \cite{ury2}.
More generally, spectral risk measures with monotonically nondecreasing spectral functions are coherent \cite{acerbi},
and on a nonatomic probability space, law invariant coherent risk measures can be represented as the supremum of a set of coherent spectral risk measures (e.g., \cite[Section 6.3.4]{SDR}).
Coherent risk measures also possess the dual representation,
which relates the risk averse and distributionally robust approaches to stochastic optimization (e.g., \cite[Section 6.3]{SDR}).
On the other hand, the Value-at-Risk, which is widely used in applications,
does not satisfy the convexity axiom and does not have such dual representation.

The convexity of risk measures, however, comes at a price from the statistical point of view.
It is well known that the Value-at-Risk is much less sensitive to outliers than the expected value (cf., \cite{huber,rous87}).
That is, when the distribution is contaminated by a small mass placed far away,
the Value-at-Risk is not changed much, while the expected value can be moved by an arbitrary amount.
It appears  that the same property of statistical nonrobustness  is shared by a large class of coherent risk measures.
Indeed, the Average Value-at-Risk and spectral coherent risk measures are concave with respect to the probability measure,
and consequently an arbitrarily small contamination of the distribution can increase their values by an arbitrary amount. Actually some authors argue that this is an advantage of the Average Value-at-Risk compared with the Value-at-Risk, since it accounts for possible large  losses (cf., \cite{ury2}). 
In any case, law invariant coherent risk measures are not robust from the statistical point of view,
which motivates us to consider risk measures that do not satisfy the convexity axiom.

For sequential decision making, an approach to risk averse multistage stochastic programming based on nested compositions of conditional risk measures was suggested in \cite{ruzshap2},
and was extended to Markov Decision Processes (MDPs) in \cite{rusz2010}.
A construction of nested risk functionals which does not rely on the convexity of risk measures was recently discussed in \cite{ShapLi},
for both Stochastic Optimal Control (SOC) and MDPs.

Closely related to the nested formulation is the time consistency of policies,
that is, whether a policy optimal at the initial stage remains optimal at every later stage conditional on the history of the process.
It is well known that  the nested formulation is time consistent,
while this is not the case for static formulations based on the Average Value-at-Risk.
It should be noted that optimality at every later stage is measured by conditional risk measures,
and hence time consistency depends on the choice of these conditional risk measures (e.g., \cite{PLS2022}).
For a given static risk measure, a natural choice is the one for which the static risk measure is equal to the composition of the conditional risk measures,
since in that case the static problem can be written in the nested form.
However,  such a choice with law invariant coherent conditional risk measures exists only for the expectation and the max-risk measure (cf.,  \cite{Kupper2009}, \cite{Shap2012}), 
and in particular does not exist for the Average Value-at-Risk.
In addition, even when such a choice exists, an optimal policy is not necessarily time consistent if the conditional risk measures are not strictly monotone \cite{ShapUgu2016,ShapXin2020}.
In this manuscript, for a given policy such that no policy is uniformly better,
we construct conditional risk measures with respect to which this policy is time consistent,
where the static risk measure is not necessarily convex.

The main contribution of this manuscript can be summarized as follows.
\begin{itemize}
  \item
  We show that risk measures concave in the probability measure, as well as their suprema, are not robust with respect to contamination of the probability measure.
  In particular, this applies to spectral coherent risk measures, such as the Average Value-at-Risk, and to law invariant coherent risk measures on nonatomic probability spaces.
  This contrasts with the Value-at-Risk, which is robust against such contamination.
  \item
  For an arbitrary law invariant risk measure satisfying the axioms of monotonicity, translation equivariance and positive homogeneity,
  and a policy such that no policy is uniformly better,
  we construct conditional risk measures with respect to which this policy is time consistent.
  The construction is based on the tail events on which the given policy is not worse than any other policy,
  and is carried out for both MDPs and SOC.
\end{itemize}

The remainder of this manuscript is organized as follows.
In Section \ref{sec-riskfun}, we discuss risk functionals which are not necessarily convex,
and their concavity and statistical robustness with respect to the probability measure.
In Sections \ref{sec_tc_mdp} and \ref{sec_tc_soc}, we construct conditional risk measures with respect to which a given policy is time consistent,
in the MDP and SOC settings, respectively.
We conclude the manuscript in Section \ref{sec_conclusion}.

We use the following notation and terminology. For $a\in \bbr$, we denote $[a]_+:=\max\{0,a\}$.
A set $\O$ equipped with its sigma algebra $\F$ is called a measurable space.
 By $\cP$ we denote  the set of all probability measures on $(\O,\F)$.
For $P\in \cP$, $(\O,\F,P)$ is called the probability space. It is said  that $P\in \cP$ is nonatomic, if for any $A\in \F$ such that $P(A)>0$, there exists $A'\in \F$ such that $A'\subset A$ and $P(A)>P(A')>0$.
By $\delta_a$ we denote the Dirac measure of mass one at the point $a$.

For $P\in \cP$ and  measurable  variable    $Z:\O\to\bbr$   we denote by $\bbe^P[Z]=\int_\O ZdP$ its expected value, and by
$$F_Z^P(z):=P(Z\le z),\;z\in \bbr,$$ its cumulative distribution function (cdf). 
By 
\[
(F_Z^P)^{-1}(\alpha):=\inf\{z:F_Z^P(z)\ge \alpha\},\;\alpha\in [0,1], 
\]
we denote the left side quantile. In particular if $Z\sim P$ we   write $F^P(z)$ and $(F^P)^{-1}(\alpha)$ for the respective cdf and quantile. 
By  $\cF$ we denote the set of right side continuous monotonically nondecreasing functions $\phi:\bbr\to \bbr$ such that $\lim_{z\to-\infty}\phi(z)=0$ and $\lim_{z\to +\infty}\phi(z)=1$, i.e. $\cF$ is the space of cdfs.
By $L_p(\O,\F,P)$ we denote the space of random variables $Z$ such that $\int_\O|Z|^p dP<\infty$, $p\in [1,\infty)$.

\setcounter{equation}{0}
\section{Risk functionals}
\label{sec-riskfun}

Let $(\O,\F)$ be a measurable space. Assume that with every probability measure  $P\in \cP$  is associated a linear space $\Z$  of measurable functions (random variables) $Z:\O\to \bbr$ and a functional $\R:\Z\to \bbr$.  Consider the following axioms (conditions) that $\R$  may  satisfy.
\begin{itemize}
  \item [(A1)] (monotonicity) if $Z,Z'\in \Z$ and $Z\ge Z'$, $P$-a.s., then $\R(Z)\ge \R(Z')$.
    \item [(A2)] (translation equivariance) for   $Z\in \Z$ and $\tau\in \bbr$, it follows that $\R(Z+\tau)=\R(Z)+\tau$.
     \item [(A3)] (positive homogeneity) if $Z\in \Z$ and $\tau\ge 0$, then
     $\R(\tau Z)=\tau \R(Z)$.
      \item [(A4)] (convexity)    if $Z,Z'\in \Z$ and $\tau\in [0,1]$, then
    $\R(\tau Z+(1-\tau)Z')\le \tau \R(Z)+(1-\tau)\R(Z')$.
      \end{itemize}
In order to emphasize that the functional $\R$ depends on $P\in \cP$,  we write it as $\R^P$. The corresponding linear space $\Z$ may also depend on $P$, but we suppress this in the notation. 
 In any case, we assume that if $Z\in \Z$ and $\tau\in \bbr$, then $Z+\tau\in \Z$. A typical example of the space $\Z$ is $L_p(\O,\F,P)$, $p\in [1,\infty]$. 
Note that it follows from axiom (A1) that if $Z=Z'$, $P$-a.s., then $\R^P(Z)= \R^P(Z')$.
It follows from axiom (A3) that $\R^P(0)=0$.

\begin{definition}
\label{def-risk}
It is said that functional $\R^P$ is a {\em risk measure} if it satisfies axioms (A1) and (A2). 
\end{definition}

Functionals satisfying axioms (A1)-(A4) are called {\em coherent} risk measures \cite{ADEH:1999}.
 An important example of a non-convex risk measure   is the Value-at-Risk, $\R^P=\VaR^P_{\alpha}$, where
\begin{equation}\label{quan}
\VaR^P_{\alpha}(Z):=\inf\{z:F_Z^P(z)\ge 1-\alpha\}=(F_Z^P)^{-1}(1-\alpha),\;\;\alpha\in (0,1).
\end{equation}
That is, $\VaR^P_{\alpha}(Z)$ is the left side $(1-\alpha)$-quantile of the distribution of $Z$. In that example the corresponding space $\Z$ consists of all measurable $Z:\O\to \bbr$.
The $\VaR^P_{\alpha}(\cdot)$ functional satisfies axioms (A1)-(A3), but not the convexity  axiom (A4).
An important example of a coherent risk measure is the Average Value-at-Risk (also called Conditional  Value-at-Risk, expected shortfall,   expected tail loss), for  $\alpha\in (0,1]$,
\begin{equation}\label{avar}
\avr^P_{\alpha}(Z)=\inf_{\tau\in \bbr}\{\tau+\alpha^{-1}\bbe^P[Z-\tau]_+\} 
 =\alpha^{-1}\int_0^{\alpha}\VaR^P_{\tau}(Z)d\tau.
\end{equation}
The first of the above representations of the Average Value-at-Risk is due to \cite{ury2}.
In that example, the respective  space is $\Z=L_1(\O,\F,P)$.

 It is said that $Z,Z'\in \Z$ are distributionally equivalent (with respect to $P\in \cP$), if their cdfs do coincide, i.e.,
 $P(Z\le z)=P(Z'\le z)$ for all $z\in \bbr$.
\begin{definition}
It is said that $\R^P$ is {\rm law invariant}   if for any distributionally equivalent  $Z,Z'\in \Z$ it follows that  $\R^P(Z)=\R^P(Z')$.
\end{definition}

That is,  a law invariant risk measure $\R^P(Z)$ is a function of the cdf $F_Z^P$, i.e.,
it can be represented as
\begin{equation}\label{cdfr}
 \R^P(Z)=\rho (F_Z^{P}),\;\;Z\in \Z,
 \end{equation}
 where $\rho$ is a mapping from the space of cdfs $\cF$ into $\bbr$.
For example, for $\R^P=\avr^P_\alpha$ the corresponding
\[
\rho(\phi)=\inf_{\tau\in \bbr}\left\{\tau+\alpha^{-1} \int_{-\infty}^{+\infty} [z-\tau]_+d \phi(z)\right\},
\]
defined for such $\phi \in \cF$ that $\int_{-\infty}^{+\infty} [z-\tau]_+d \phi(z)<\infty$ for all $\tau\in \bbr$.
For $\R^P=\VaR_\alpha^P$ the corresponding mapping is
$$\rho(\phi)=\inf\{z:\phi(z)\ge 1-\alpha\},\;
\phi\in \cF.$$


Risk measure $\R$  is called {\em spectral} if it can be represented in the form 
\begin{equation}\label{spect}
 \R^P(Z)=\int_0^1 \sigma(\tau) (F_Z^P)^{-1}(\tau)d\tau,
 \end{equation}
where  $\sigma(\cdot)$ is a density function on [0,1]. Spectral risk measures were introduced in Acerbi \cite{acerbi}. Spectral risk measures are  law invariant. It was shown in \cite{acerbi} that  the spectral risk measure is coherent in the sense of \cite{ADEH:1999}, if  the corresponding density $\sigma(\cdot)$ is monotonically nondecreasing  on  [0,1]. 
  When the space $\Z=L_p(\O,\F,P)$, it is assumed that the spectral function  $\sigma$ belongs to the dual space  $\Z^*=L_q(\O,\F,P)$.

\begin{remark}
\label{rmk_kusuoka}
Assuming that the probability space $(\O,\F,P)$ is nonatomic, it is possible to show that any law invariant coherent risk measure   can be represented as the supremum of a set of  coherent  spectral risk measures  (e.g., \cite[Section 6.3.4]{SDR}).
This motivates us to consider risk measures of the form
\begin{equation}
  \label{riskmes}
  \R^P(Z):=\sup_{\mu \in \cM} \int_0^1  (F_Z^P)^{-1}(\tau)d\mu(\tau),
\end{equation}
where $\cM$ is a set of probability measures on the interval [0,1].
\end{remark}

\begin{definition}
\label{def-genrisk}
   We refer to 
risk measures of the form \eqref{riskmes} as the {\em generalized} risk measures. We say that a generalized risk measure is {\em generalized spectral} if the set $\cM=\{\bar{\mu}\}$ is a singleton.  
\end{definition}

The generalized risk measure also depends on a choice of the set $\cM$. In order to emphasize this, we sometimes write it as $\R^P_\cM(\cdot)$.
The generalized risk measures are law invariant and 
satisfy axioms (A1)-(A3), but may not satisfy the convexity   axiom (A4).  
 By the above we have that a generalized  risk measure is coherent if the probability measures $\mu\in \cM$ have monotonically nondecreasing  densities.   An important example of a non-convex generalized spectral risk measure is 
$\VaR^P_{\alpha}$, which  can be written in the form  \eqref{riskmes} with $\cM$ being the singleton $\cM=\{\delta_{1-\alpha}\}$.

\subsection{Concavity of risk measures}

\begin{definition}
\label{def-concave}
A law invariant functional $\R^{P}:\Z\to\bbr$ is said to be  concave in $P$,     if for  $P, P'\in \cP$ and $\tau\in [0,1]$ it follows   for all $Z\in \Z$ that 
 \begin{equation}\label{concavefunc}
 \R^{\tau P+(1-\tau) P'}(Z )\ge \tau \R^P(Z)+ (1-\tau)\R^{P'}(Z).
 \end{equation}
\end{definition}

Such concavity property of risk measures is useful, for example, in construction of upper bounds in risk averse SDDP type algorithms \cite{upperbound}, and in establishing existence of non-randomized policies in risk averse MDPs \cite{ShapLi}.

 There is a large class of   risk functionals which are concave with respect to the probability distribution. For example consider   risk functionals of the form $\R^P(Z)=\inf_{\theta\in \Theta}\bbe^P[\Psi(Z,\theta)]$, where $\Theta$  is a subset of a finite dimensional vector space and
$\Psi:\bbr\times \Theta\to \bbr$ is a real valued function. The $\avr_\alpha$ risk measure is of that form.  In that case  $\R^P(Z)$ is given by the infimum of functionals linear in $P$, and hence  is concave in $P$.
 Spectral coherent risk measures are representable as the corresponding integral of the Average Value-at-Risk measures and hence are concave in $P$.

We refer to the interval [0,1] equipped with its Borel sigma algebra   as the {\em standard} sample space. We can view generalized risk measures as defined on the standard sample space. Consider the Bernoulli probability measure $B_\tau=\tau\delta_0 + (1-\tau)\delta_1$, $\tau\in [0,1]$. With a law invariant risk measure $\R^P$ and $Z\in \Z$ we can associate the function
\begin{equation}\label{bern}
 \cp_Z(\tau):=\R^{B_\tau}(Z),\;\tau\in [0,1]. 
\end{equation}
Since for $\gamma,\tau_1,\tau_2\in [0,1]$,   
\[
B_{\gamma \tau_1+(1-\gamma)\tau_2}=\gamma(\tau_1\delta_0+(1-\tau_1)\delta_1)+
(1-\gamma)(\tau_2\delta_0+(1-\tau_2)\delta_1)=\gamma B_{\tau_1}+ (1-\gamma)B_{\tau_2},
\]
it follows that if $\R^P$ is concave in $P$, then the function $\cp_Z(\cdot)$ is concave on [0,1] for all $Z\in \Z$. For    $Z\sim B_\tau$   we write  $\cp(\tau)$ for the corresponding function. 

Let us consider some examples. Suppose that $Z\sim B_\tau$. Then  we have that $F^{B_\tau}(z)=0$ for $z<0$, 
$F^{B_\tau}(z)=\tau$ for $z\in [0,1)$ and $F^{B_\tau}(z)=1$ for $z\ge 1$. Thus
$(F^{B_\tau})^{-1}(\beta)=0$ for $\beta\in [0,\tau)$, and $(F^{B_\tau})^{-1}(\beta)=1$ for $\beta\in [\tau,1]$.
 Consider   $\R^P(Z):= \VaR^P_{\alpha}(Z)=(F_Z^P)^{-1}(1-\alpha)$ for some $\alpha\in (0,1)$. Then for  $Z\sim B_\tau$,
 \begin{equation}\label{func-1}
 \cp(\tau)=\left\{
 \begin{array}{lll}
 0,\;{\rm for}\;\tau\in [0,\alpha),\\
1,\;{\rm for}\; \tau\in [\alpha,1].
 \end{array}
 \right.
 \end{equation}
 Here $\cp(\cdot)$ cannot be represented as the supremum of functions concave on the interval $[0,1]$, and hence 
   $\VaR^P_{\alpha}(Z)$ cannot be represented as the supremum of law invariant risk measures concave in $P$. 

Consider now  $\R^P(Z):= \avr^P_{\alpha}(Z)$. In that case by   \eqref{avar} we have that  for  $Z\sim B_\tau$,
\begin{equation}\label{func-2}
 \cp(\tau)=\left\{
 \begin{array}{lll}
\alpha^{-1} \tau,\;{\rm for}\;\tau\in [0,\alpha),\\
1,\;{\rm for}\;\tau\in [\alpha,1].
 \end{array}
 \right.
 \end{equation}
 Of course this function is concave on $[0,1]$.

\subsection{Statistical robustness of risk functionals }\label{sec_stat_rob}

It is well known that  the Value-at-Risk (quantile) functionals $\VaR^P_{\alpha}$ are much less sensitive to outliers than the expected value $\bbe^P$. That is, for a probability measure $P$ supported  on the  real line and   $\e\in (0,1)$ consider the contamination $(1-\e)P+\e \delta_\kappa$ of the original distribution $P$ by the scaled Dirac measure. Then   
\[
\bbe^{(1-\e)P+\e \delta_\kappa}[\cdot]=(1-\e)\bbe^P[\cdot]+\e \kappa.
\]
By choosing $\kappa\in \bbr$  large enough, the contaminated expectation can be moved from the expectation $\bbe^P$ by an arbitrary amount for any positive value of $\e$.   On the other hand the $\VaR^P_{\alpha}$ measure  will not be changed much by contamination for $\e<\alpha$. The respective theory of the so-called Robust Statistics was well discussed already more than  40 years ago (cf., \cite{huber},\cite{rous87}). 

Suppose now that risk measure $\R^P$ satisfies axioms (A1)-(A3) and is concave in $P$. Then for $\e\in (0,1)$, 
\begin{equation}\label{inequal}
\R^{(1-\e)P+\e \delta_\kappa}(\cdot)\ge (1-\e) \R^P(\cdot)+ \e \R^{\delta_\kappa}(\cdot) = (1-\e)\R^P(\cdot)+\e \kappa,
\end{equation}
where the last equality holds since $\delta_\kappa$ is the distribution of the constant $\kappa$, and hence $\R^{\delta_\kappa}(\cdot)=\kappa$ by axioms (A2)-(A3).
As in the expected value example,   the risk  $\R^{(1-\e)P + \e \delta_\kappa} (\cdot)$ can be moved from $\R^P(\cdot)$ by an arbitrary amount for any positive value of $\e$. As was already mentioned, spectral coherent risk measures are concave in $P$, and hence are not robust from the statistical point of view. Moreover, consider functionals $\cR^P(\cdot):=\sup_{i\in I}\R_i^P(\cdot)$, where $\{\R_i^P\}_{i\in I}$ is a family of risk measures such that for every $i\in I$ the risk measure $\R_i^P(\cdot)$ is law invariant, satisfies axioms (A1)-(A3) and is concave in $P$.
Then   by \eqref{inequal}
\[
\begin{array}{lll}
\cR^{(1-\e)P+\e \delta_\kappa}(\cdot) = \sup\limits_{i\in I}\R_i^{(1-\e)P+\e \delta_\kappa}(\cdot) 
\ge  \sup\limits_{i\in I}\{(1-\e)\R_i^P(\cdot)+  \e \kappa\} 
= (1-\e)\cR^P(\cdot)+\e \kappa,
\end{array}
\]
and hence $\cR^P$ is not robust from the statistical point of view.
In particular, it follows from Remark \ref{rmk_kusuoka} that law invariant coherent risk measures on a nonatomic probability space are not robust from the statistical point of view.
 
\section{Time consistency}

In this section, we discuss time consistency of policies with respect to nested formulations of risk averse problems.
Specifically, in Section \ref{sec_tc_mdp}, for a given policy of an MDP such that no policy is uniformly better,
we construct conditional risk measures with respect to which this policy is time consistent,
and in Section \ref{sec_tc_soc}, we carry out a similar construction for SOC.

\subsection{Markov Decision Processes}\label{sec_tc_mdp}
 
Consider a finite horizon MDP with  state spaces $\s_t$,  action spaces $\A_t$, cost functions $c_t(s_t,a_t,s_{t+1})$, and transition probabilities $P_t(s_{t+1}|a_t,s_t)$, 
  $t=1,...,T$. A (non-randomized) policy $\pi$ is defined by mappings $\pi_t:(s_1,...,s_t)\mapsto \A_t$. A policy is Markovian if $\pi_t(s_t)$, $t=1,...,T$,  are functions of $s_t$ rather than the whole history $(s_1,...,s_t)$.  
To a policy $\pi$ corresponds the  probability distribution    defined on histories $h_t=(s_1,...,s_t)$ with $a_\tau=\pi_\tau(s_1,...,s_\tau)$, $\tau=1,...,t$ (cf., \cite{IT1949}).

 We say  that a policy $\bar{\pi}$ is {\em time consistent} if at every stage $t$, conditional on the history $h_t$, the policy $\{\bar{\pi}_t,...,\bar{\pi}_T\}$   is optimal. 
 In order to make this precise we need to define what we mean by optimality at every stage $t$. Suppose that for every $t$ the state space $\s_t$ is a measurable space. 
  Consider the set $\Pi$ of   Markovian policies. 
Let  $\R^P$ be a law invariant risk measure with $P$ being a probability measure defined on $(s_1,\dots, s_T)\in \s_1\times\cdots\times \s_T$. 
Let $\bbp_t^{s_t,a_t}$, $t=1,...,T$,  be reference transition kernels. Then to a policy $\pi\in \Pi$ corresponds the nested functional $\cR^\pi$ given by the composition of conditional risk measures $\R^{\bbp_t^{s_t,\pi_t(s_t)}}$ (cf., \cite[Section 4]{ShapLi}). 
Let us make the following observation.
\begin{remark}
\label{rem-a}
  Let $\cP$ be the set of probability measures on a measurable space  $(\O,\F)$  and $\R^P$ be a law invariant risk measure associated with $P\in \cP$. For a given $\bbp\in \cP$ and a set $C\in \F$ such that $\bbp(C)>0$, consider the induced measure $\sP$ supported on $C$, i.e., 
$\sP(A)=\bbp(A)/\bbp(C)$ for any measurable set $A\subset C$. Consider the corresponding risk measure $\R^{\sP}$. Then for any measurable $Z,Z':\O\to \bbr$ such that $Z(\w)=Z'(\w)$ for all $\w\in C$, we have  that $Z$ and $Z'$ are distributionally equivalent with respect to $\sP$ and hence   $\R^{\sP}(Z)=\R^{\sP}(Z')$.
\end{remark}

This motivates the following construction.
  Let $\bar{\pi}\in \Pi$ be a given policy (for example, $\bar{\pi}$ can be an optimal policy with respect to $\cR^\pi$ for the total cost). 
 Assume  that  there is no policy uniformly better than $\bar{\pi}$, that is: 
\begin{assumption}
\label{ass-1}
For every policy   $\{\pi_t\}$ and all  $s_{t}\in \s_{t}$ there is a set $C_t=C_t(\pi_t,s_{t})\subset \s_{t+1}$ of positive measure with respect to the  reference probability measure $\bbp_t^{{s_t,\pi_t(s_t)}}$ on $\s_{t+1}$,   such that $Z^{\pi_t}(s_{t+1})\ge \bar{Z} (s_{t+1})$ for a.e. $s_{t+1}\in C_t$, where $$Z^{\pi_t}(s_{t+1}):= c_t(s_t,\pi_t(s_t),s_{t+1}) \;\;{\rm and}\;\;\bar{Z}(s_{t+1}):= c_t(s_t,\bar{\pi}_t(s_t),s_{t+1}).$$
\end{assumption}

Consider the probability measure $\sP_{t}$ on $C_t$  induced by $\bbp_t^{{s_t,\pi_t(s_t)}}$. That is, for a measurable set $A\subset C_t$, $$\sP_{t}(A)=\bbp_t^{{s_t,\pi_t(s_t)}}(A)/\bbp_t^{{s_t,\pi_t(s_t)}}(C_t).$$   
Note that the measure $\sP_{t}=\sP_{t}^{{s_t,\pi_t(s_t)}}$ is supported on the set $C_t=C_t(\pi_t,s_{t})$ and depends on $\pi_t$ and $s_{t}$.  Consider the corresponding conditional risk measures $\R^{\sP_{t}^{s_t,\pi_t(s_t)}}$.
In the following, we show that under Assumption \ref{ass-1}, policy $\bar{\pi}$ is time consistent with respect to the nested formulation determined by risk measures constructed from $\R^{\sP_{t}^{{s_t,\pi_t(s_t)}}}$.

Specifically, we construct risk measures $\widetilde{\R}_t^{s_t}$ such that $\bar{\pi}$ is time  consistent with respect to the corresponding nested formulation. 
Throughout the rest of our discussion, we assume that $\R$ also satisfies axiom (A3), and hence $\R^{\sP}(0) = 0$ for any probability measure $\sP$.
With the measures $\sP_t^{s_t, \pi_t(s_t)}$ and the function $\bar{Z}$ defined above,
we define the risk measure $\widetilde{\R}_t^{s_t}$ of $Z: \s_{t+1} \to \bbr$ as
\begin{align}
\label{tail_risk_mdp}
\widetilde{\R}_t^{s_t}(Z)
:= \sup_{\pi \in \Pi} \R^{\sP_t^{s_t, \pi_t(s_t)}} \big( Z - \bar{Z} \big).
\end{align}
For any $\pi \in \Pi$, we can write the value function of $\pi$ under the nested formulation determined by $\widetilde{\R}_t^{s_t}$ recursively as
$V^\pi_{T+1}(\cdot) \equiv 0$ and
\begin{align}
\label{value_pi_mdp}
V^\pi_t(s_t)
= \widetilde{\R}_t^{s_t} \big( c_t(s_t,\pi_t(s_t),\cdot) + V^\pi_{t+1}(\cdot) \big),
~ s_t\in \s_t, ~ t=T,\ldots,1.
\end{align}

\begin{proposition}
\label{prop_tc_mdp}
Suppose Assumption \ref{ass-1} holds. 
Then for any $\pi \in \Pi$, $t=1,\ldots,T$ and $s_t\in \s_t$, we have
\begin{align*}
V^{\bar{\pi}}_t(s_t) \leq V^{\pi}_t(s_t).
\end{align*}
That is, $\bar{\pi}$ is time consistent with respect to the nested formulation determined by $\widetilde{\R}_t^{s_t}$.
\end{proposition}

\begin{proof}
Since $V^{\bar{\pi}}_{T+1} \equiv 0$ and $\widetilde{\R}_t^{s_t}(\bar{Z}) = 0$,
we obtain from \eqref{value_pi_mdp} with $\pi=\bar{\pi}$ that $V^{\bar{\pi}}_t \equiv 0$ for $t = T, \ldots, 1$.
Clearly the desired claim holds at $t = T+1$. 
Suppose the claim holds at $t+1$.  Then we have 
\begin{align*}
V^{\bar{\pi}}_t(s_t)
= \R^{\sP_t^{s_t, \pi_t(s_t)}}(0)
& \overset{(a)}{\leq} 
\R^{\sP_t^{s_t,\pi_t(s_t)}} \big(Z^{\pi_t} - \bar{Z}\big) \\
& \overset{(b)}{\leq}
\widetilde{\R}_t^{s_t} \big( Z^{\pi_t} \big) \\
& \overset{(c)}{=}
\widetilde{\R}_t^{s_t} \big( c_t(s_t, \pi_t(s_t), \cdot) + V^{\bar{\pi}}_{t+1}(\cdot) \big) \\
& \overset{(d)}{\leq}
\widetilde{\R}_t^{s_t} \big( c_t(s_t, \pi_t(s_t), \cdot) + V^\pi_{t+1}(\cdot) \big)
= V^\pi_t(s_t),
\end{align*}
where (a) follows from Assumption \ref{ass-1} and (A1), (b) follows from \eqref{tail_risk_mdp},
(c) follows from $V^{\bar{\pi}}_{t+1} \equiv 0$, and (d) follows from the induction hypothesis.
Hence the claim holds at $t$.
The proof is then completed.
\end{proof}

It could be noted that the construction of $\widetilde{\R}_t^{s_t}$ does not rely on the policies being Markovian.
Consider history dependent policies $\pi_t(h_t)$ with $h_t = (s_1, \ldots, s_t)$,
and suppose that for every policy $\pi$ and every history $h_t$ there is a set $C_t = C_t(\pi_t, h_t) \subset \s_{t+1}$
of positive measure with respect to $\bbp_t^{s_t, \pi_t(h_t)}$
such that $c_t(s_t, \bar{\pi}_t(h_t), s_{t+1}) \leq c_t(s_t, \pi_t(h_t), s_{t+1})$ for a.e. $s_{t+1} \in C_t$.
Then the conclusion of Proposition \ref{prop_tc_mdp} holds for history dependent policies,
by following the same construction with $s_t$ replaced by $h_t$.
The same applies to the SOC setting discussed in Section \ref{sec_tc_soc} below.

\subsection{Stochastic Optimal Control}\label{sec_tc_soc}

Let us now consider the  stochastic optimal control (SOC) problem 
\begin{align*}
\min_{\pi \in \Pi} \bbe^\pi \Big[ \sum_{t=1}^{T} c_t(x_t, u_t, \xi_t) \Big],
\end{align*}
where the state evolves according to $x_{t+1} = \Phi_t(x_t, u_t, \xi_t)$, $t = 1, \ldots, T$, from a deterministic initial state $x_1$,
and the minimization is over policies that determine the controls $u_t\in \U_t$. 
Here $\X_t \subset \bbr^{n_t}$, $\U_t \subset \bbr^{m_t}$ and $\Xi_t \subset \bbr^{d_t}$ are nonempty closed sets,
$\xi_t \in \Xi_t$ is the random data at stage $t$, $\F_t$ is the Borel sigma algebra of $\Xi_t$,
and both the cost function $c_t: \X_t \times \U_t \times \Xi_t \to \bbr$ and the mapping $\Phi_t: \X_t \times \U_t \times \Xi_t \to \X_{t+1}$ are continuous.

Throughout the rest of our discussion, we assume that the distribution of $(\xi_1, \ldots, \xi_T)$ does not depend on states and controls,
and is given by the product $P = P_1 \times \cdots \times P_T$ of the marginal distributions $P_t$ of $\xi_t$ on $(\Xi_t, \F_t)$.
For each $P_t$,  let us consider a law invariant risk measure $\R^{P_t}$ satisfying axioms (A1)-(A3).
The value functions of the risk averse counterpart of the SOC problem then depend only on the state, and satisfy the dynamic programming equations
$V_{T+1}(\cdot) \equiv 0$ and
\begin{align*}
V_t(x_t) = \inf_{u_t \in \U_t} \R^{P_t} \big[ c_t(x_t, u_t, \xi_t) + V_{t+1} \big( \Phi_t(x_t, u_t, \xi_t) \big) \big],
~ t = T, \ldots, 1,
\end{align*}
where $\xi_t \sim  P_t$.
In view of this, we restrict our attention to the set $\Pi$ of Markovian policies $u_t = \pi_t(x_t)$, $t = 1, \ldots, T$.

Similar to Section \ref{sec_tc_mdp}, for a given policy $\bar{\pi} \in \Pi$,
we proceed to construct risk measures $\widetilde{\R}_t^{x_t}$ such that $\bar{\pi}$ is time consistent with respect to the corresponding nested formulation.

\begin{assumption}
\label{ass_tail_soc}
For $t=T,\ldots,1$, any policy $\pi\in\Pi$ and any $x_t \in \X_t$, there is a set $C_t = C_t(\pi_t, x_t) \subset \Xi_t$ such that
$P_t(C_t) > 0$ and $\bar{Z}(\xi_t) \leq Z^{\pi_t}(\xi_t)$ for a.e. $\xi_t \in C_t$, where
\begin{align*}
Z^{\pi_t}(\xi_t) := c_t(x_t,\pi_t(x_t),\xi_t)
~~ {\rm and} ~~
\bar{Z}(\xi_t) := c_t(x_t, \bar\pi_t(x_t),\xi_t).
\end{align*}
\end{assumption}

For $\pi \in \Pi$ and $x_t \in \X_t$, consider the probability measure $\sP_{C_t}$ on $C_t = C_t(\pi_t, x_t)$ induced by $P_t$,
that is, $\sP_{C_t}(A)=P_t(A)/P_t(C_t)$ for measurable $A \subset C_t$.
We then define the risk measure $\widetilde{\R}_t^{x_t}$ of $Z: \Xi_t \to \bbr$ as
\begin{align}
\label{tail_risk_soc}
\widetilde{\R}_t^{x_t}(Z)
:= \sup_{\pi \in \Pi} \R^{\sP_{C_t(\pi_t, x_t)}} \big( Z - \bar{Z} \big).
\end{align}
For any $\pi \in \Pi$, we can write the value function of $\pi$ under the nested formulation determined by $\widetilde{\R}_t^{x_t}$ recursively as
$V^\pi_{T+1}(\cdot) \equiv 0$ and
\begin{align}
\label{value_pi_soc}
V^\pi_t(x_t)
= \widetilde{\R}_t^{x_t} \Big( c_t(x_t, \pi_t(x_t), \cdot) + V^\pi_{t+1} \big( \Phi_t(x_t, \pi_t(x_t), \cdot) \big) \Big),
~ x_t \in \X_t, ~ t = T, \ldots, 1.
\end{align}

\begin{proposition}
\label{prop_tc_soc} 
Suppose Assumption \ref{ass_tail_soc} holds.
Then for any $\pi \in \Pi$, $t = 1, \ldots, T$ and $x_t \in \X_t$, we have
\begin{align*}
V^{\bar{\pi}}_t(x_t) \leq V^{\pi}_t(x_t).
\end{align*}
That is, $\bar{\pi}$ is time consistent with respect to the nested formulation determined by $\widetilde{\R}_t^{x_t}$.
\end{proposition}

\begin{proof}
Similar to the proof of Proposition \ref{prop_tc_mdp},
we obtain from \eqref{value_pi_soc} with $\pi = \bar{\pi}$ and \eqref{tail_risk_soc} that $V^{\bar{\pi}}_t \equiv 0$ for $t = T, \ldots, 1$.
Clearly the desired claim holds at $t = T+1$.
Suppose the claim holds at $t+1$. With $C_t = C_t(\pi_t, x_t)$, we have
\begin{align*}
V^{\bar{\pi}}_t(x_t)
= \R^{\sP_{C_t}}(0)
& \overset{(a)}{\leq}
\R^{\sP_{C_t}}\big( Z^{\pi_t}-\bar{Z} \big) \\
& \overset{(b)}{\leq}
\widetilde{\R}_t^{x_t} \big( Z^{\pi_t} \big) \\
& \overset{(c)}{=}
\widetilde{\R}_t^{x_t} \Big( c_t(x_t, \pi_t(x_t), \cdot) + V^{\bar{\pi}}_{t+1} \big( \Phi_t(x_t, \pi_t(x_t), \cdot) \big) \Big) \\
& \overset{(d)}{\leq}
\widetilde{\R}_t^{x_t} \Big( c_t(x_t, \pi_t(x_t), \cdot) + V^\pi_{t+1} \big( \Phi_t(x_t, \pi_t(x_t), \cdot) \big) \Big)
= V^\pi_t(x_t),
\end{align*}
 where (a) follows from Assumption \ref{ass_tail_soc} and (A1), (b) follows from \eqref{tail_risk_soc}, 
(c) follows from $V^{\bar{\pi}}_{t+1} \equiv 0$,
and (d) follows from the induction hypothesis with $x_{t+1} = \Phi_t(x_t, \pi_t(x_t), \xi_t)$ for any $\xi_t \in \Xi_t$.
Hence the claim holds at $t$.
The proof is then completed.
\end{proof}

\section{Concluding remarks}\label{sec_conclusion}

In this manuscript, we discuss risk measures which are not necessarily convex.
We show that concavity in the probability measure, which holds for a large class of coherent risk measures,
makes the risk sensitive to an arbitrarily small contamination of the distribution,
and hence such risk measures are not robust from the statistical point of view.
Our findings suggest that when statistical robustness is of concern, it is natural to consider risk measures without the convexity axiom, such as the Value-at-Risk.

For any law invariant risk measure satisfying axioms (A1)-(A3) and any policy such that no policy is uniformly better,
we construct conditional risk measures with respect to which this policy is time consistent,
in both the MDP and SOC frameworks.
It should be noted that the constructed risk measures depend on the given policy and are not law invariant.
The important question of how to solve the resulting non-convex problems numerically is beyond the scope of this manuscript,
and can be a topic of future research.

\paragraph{Acknowledgments}
We used Claude (Anthropic) to produce initial drafts of selected paragraphs,
based on the main results provided by the authors,
and to polish the writing and correct typos.
All content was subsequently revised and verified by the authors.

\bibliographystyle{plain}
\bibliography{references}

\end{document}